# On the number of mechanical configurations for nonlinear stiffness systems designed based on a linear spring with a nonlinear boundary

Qiangqiang Li[*], Qingjie Cao

*School of Astronautics, Harbin Institute of Technology, Harbin 150001, China*

[*]Corresponding author

*E-mail address:* <u>18slqq@alu.hit.edu.cn</u>

**Abstract**

In this paper, we present a novel integrated method for designing nonlinear stiffness systems based on a general spring-boundary model (GSBM) to study the number of mechanical configurations for nonlinear stiffness systems designed by the combination of a linear spring with a nonlinear boundary. GSBM consists of a lumped mass, a special-shaped track, a roller rolling in the track and GLSM with either positive or negative stiffness. The integrated method considers pre-tensioned, pre-compressed and original length conditions of GLSM to design roller trajectories to customize nonlinear stiffness systems. It is proved that the mechanical configurations of nonlinear stiffness systems designed by the combination of a linear spring with a nonlinear boundary are not limited to one, but six or eight forms: for systems with nonnegative or nonpositive potential energy, there are six independent mechanical configurations, and for other systems, there are eight independent mechanical configurations.

# 1. Introduction

Nonlinear systems exist widely in nature, however, how to construct them with accurate expected nonlinearity artificially remains a problem. In recent years, owning to that the nonlinear stiffness characteristics directly determine the dynamical performances of the nonlinear mechanical systems applied in the fields such as vibration isolation, vibration energy harvesting, nonlinear energy sink and robotic engineering, designing mechanical systems with expected nonlinear stiffness characteristics has received much attention [1-6].

Up to now, there are two main ways to design mechanical systems with expected nonlinear stiffness characteristics. In the first method, a linear spring is connected with a string wrapping around a non-circular pulley which rotates linearly, leading to a nonlinear torque [2, 3]. By designing the pulley profile, the expected nonlinear torque can be obtained. This work can date back to Michel Jean, a French scientist, who firstly used this method to construct an accurate mechanical structure of Duffing system in 1956 [7]. In the second method, a linear spring is connected with a roller which can slide along a special-shaped cam [1, 4-6]. With the cam moving in the direction perpendicular to the deformation of the linear spring, we can have a nonlinear restoring force, and it can be customized accurately by designing the shape of the cam. Collectively, the key constructing idea of these two methods is introducing a nonlinear boundary (pulley, cam, etc.) to a linear spring to produce expected nonlinear force. By contrast, the mechanism designed by the first method can only produce uni-directional force because the string can only be pulled. To produce bi-directional force, two antagonistic mechanisms must be assembled together, which will increase the complexity of the whole structure. Therefore, the second method is more feasible and effective in designing mechanical systems with expected nonlinear stiffness characteristics.

Although, using the second method, for a nonlinear stiffness system, we can always find a corresponding mechanical configuration based on the basic combination of a linear spring with a nonlinear boundary, it is still unclear how many independent configurations we can find, which is an important scientific problem that needs to be further discussed, and has not been reported yet. The motivation of this paper is to present a novel integrated second method for designing nonlinear stiffness systems based on a general spring-boundary model (GSBM), to study the number of independent mechanical configurations of nonlinear stiffness systems designed by the combination of a linear spring with a nonlinear boundary.

The rest of this paper is organized as follows. In Sec.2, a general linear spring model (GLSM) is proposed, and GSBM is established based on GLSM. In Sec. 3, the integrated method of designing nonlinear stiffness systems is developed based on GSBM. In Sec. 4, the number of mechanical configurations of nonlinear stiffness systems designed by the combination of a linear spring with a nonlinear boundary is discussed. In Sec. 5, an example of softening Duffing system is constructed. Finally, conclusions of this study are drawn in Sec. 6.

## 2. Modeling
### 2.1. GLSM

Fig. 1 shows the model of GLSM constructed by connecting a linear spring (positive stiffness element) and a dipteran flight mechanism (negative stiffness element) given by [8] in parallel. GLSM consists of a vertical spring with stiffness $K_1$, a pair of rigid rods with length $L$ and a pair of oblique springs with stiffness $K_2$ and original length $L_0$. The whole system is initially at the equilibrium position with vertical spring unstressed and the two rods with no angle about the horizontal direction, as shown in Fig. 2(a). With an external force $F_{e1}$ applied to roller $P$ in $Y$ direction, as shown in Fig. 2(b), the roller deviates from the equilibrium position with a displacement $Y$. When $Y=L$, the roller will be locked by the two rods, as illustrated in Fig. 2(c), so the effective working range of GLSM is limited within $|Y|<L$.

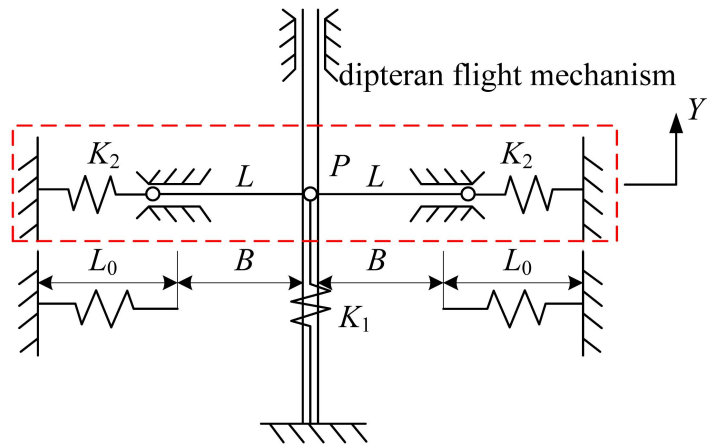

Fig. 1. Schematic diagram of GLSM.

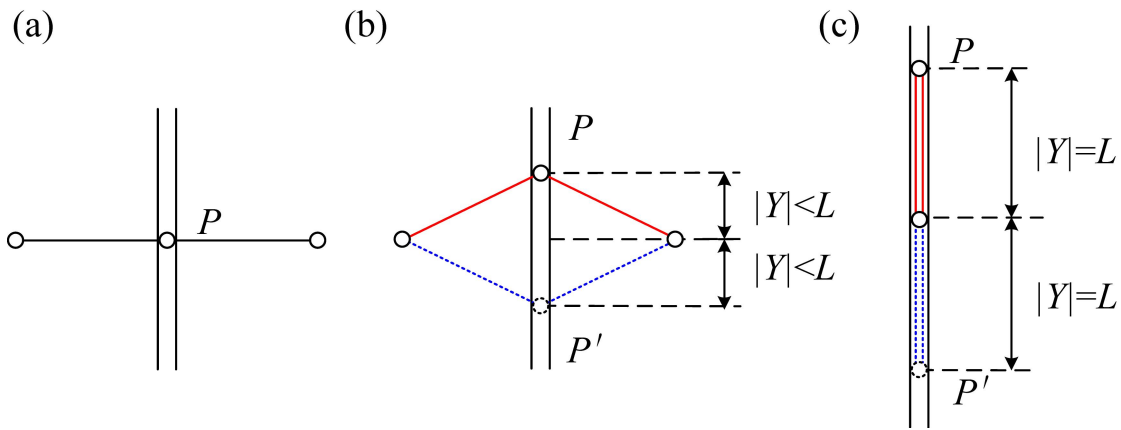

Fig. 2. Three different static states of GLSM: (a) at equilibrium position; (b) deviated from equilibrium position with displacement $Y$; (c) at 'locked' position.

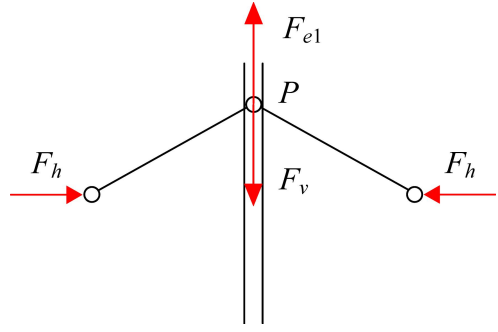

Fig. 3. Force analysis diagram of GLSM.

According to the force analysis illustrated in Fig. 3, the relationship between the applied force and the displacement can be derived as

$$F_{e1}(Y)=K_1Y-2K_2\left(1-\frac{B}{\sqrt{L^2-Y^2}}\right)Y, |Y|<L, \qquad (1)$$

where $B$ is half length between the inner ends of the horizontal springs at their un-stressed positions. By differentiating Eq. (1) with respect to displacement $Y$, the stiffness of GLSM can be given as

$$K(Y)=K_1-2K_2+\frac{2K_2BL^2}{(L^2-Y^2)^{3/2}}, |Y|<L. \qquad (2)$$

In the case when $B=0$, Eq. (2) becomes

$$K=K_1-2K_2. \qquad (3)$$

It can be seen from Eq. (3) that the stiffness of GLSM is determined by two parameters: the vertical spring stiffness $K_1$ and the oblique spring stiffness $K_2$. By adjusting the two parameters, we can have negative, zero or positive linear stiffness, which means GLSM can work as a general linear spring. In the case when GLSM works as a positive spring ($K_1>2K_2$), it is obvious that the working mode of GLSM is the same as that of the traditional spring. Specially, in the case when GLSM works as a negative spring ($K_1<2K_2$), the spring force acts in the same direction with the spring deformation, which is opposite to the working mode of the traditional spring. Note that if $K_1=2K_2$, GLSM will evolve into a zero-stiffness spring which is not the main concern in this study, thus the value range of $K$ is set to $\{K\mid K\neq 0\}$.

**2.2. GSBM**

Based on GLSM, GSBM is established, as shown in Fig. 4(b), which comprises a lumped mass, a special-shaped track (consolidated with the mass), a roller rolling along the track and GLSM. In contrast to the special-shaped-cam-spring mechanism (SSCSM) shown in Fig. 4(a), based on which the existing second method is developed, GSBM is a general model of the basic combination of a linear spring with a nonlinear boundary, lying in: 1) the track is a general boundary that can not only be compressed but also be pulled by the linear spring; 2) GLSM is a general spring with either positive or negative stiffness. In practical design, the roller trajectory of GSBM can be designed artificially to change its mechanical configuration, to customize the mechanical system with expected stiffness characteristics.

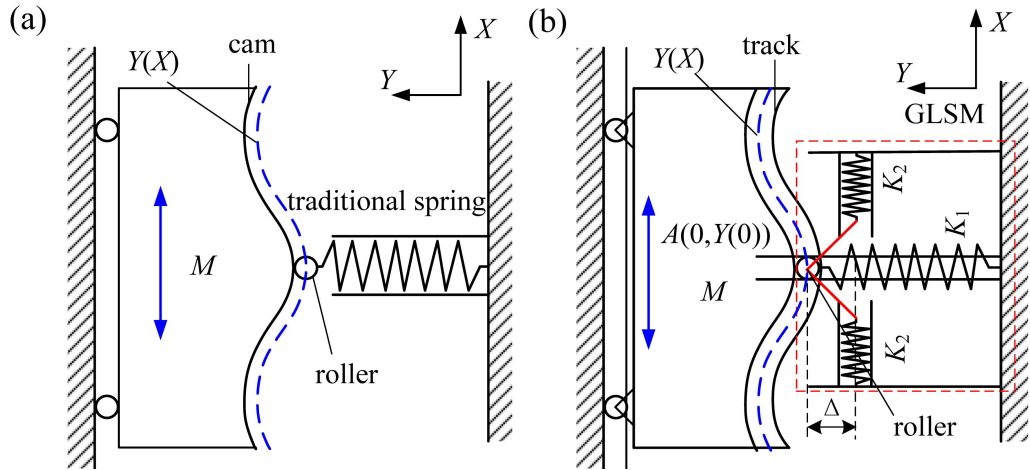

Fig. 4. Physical models of (a) SSCSM proposed in [1] and (b) GSBM.

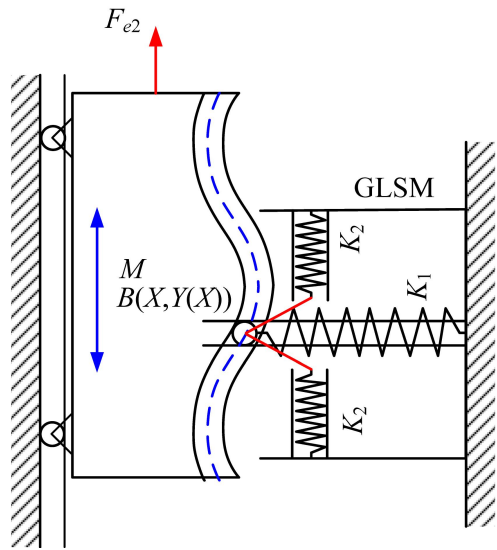

Fig. 5. GSBM deviated from the initial position with an external force.

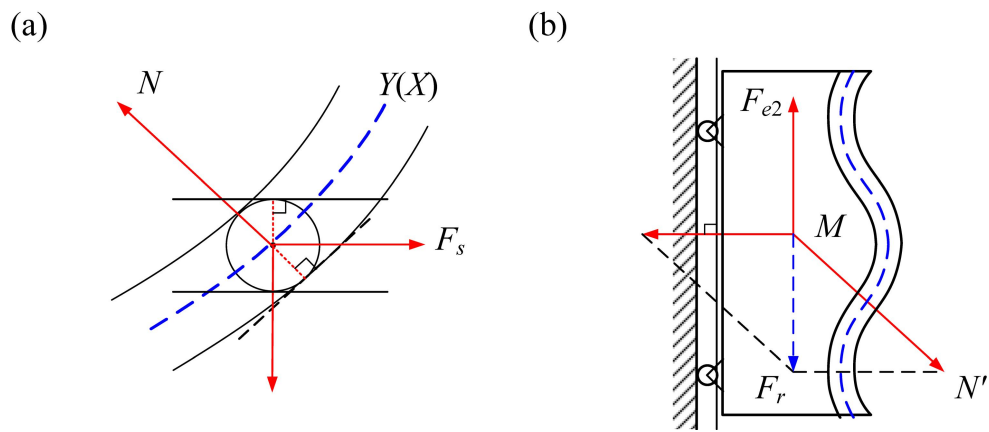

Fig. 6. Force analysis diagrams of (a) the roller and (b) the mass.

To give a clear description of the static characteristics of GSBM, an external force $F_{e2}$ is applied to the lumped mass in $X$ direction. Without the applied force, as shown in Fig. 4(b), GSBM stays at the initial position with GLSM pre-stretched (or pre-compressed) by length $\Delta$. Under the applied force, the mass moves by a displacement $X$, and the roller roles from the initial point $A$ to point $B$, which stretches (or compresses) GLSM by a displacement $Y(X)$, as illustrated in Fig. 5. In this process, though GLSM works as a linear spring, a nonlinear spring force $F_s$ with respect to displacement $X$ happens in $Y$ direction due to the nonlinear track, as illustrated by Fig. 6(a). The nonlinear track plays a vital role not only in no-linearizing GLSM but also in converting the nonlinear spring force to the motion direction by the reacting force $N'$ of supporting force $N$, as illustrated in Fig. 6(b), which results in a nonlinear restoring force $F_r$ equaling to the applied force but in opposite direction.

Based on the static analysis above, under the applied force, the nonlinear spring force of GLSM with respect displacement $X$ can be given as

$$F_s(X) = -K(K_1, K_2) Y(X). \tag{4}$$

By converting $F_s$ to the motion direction, the restoring force of GSBM can be derived as

$$F_r(X) = -K(K_1, K_2) Y(X) Y'(X), \tag{5}$$

where the prime over $Y$ denotes its derivative with respect to displacement $X$. Assuming the damping and friction of the system are not considered, the equation of motion of GSBM with no external excitation can be given as

$$M\ddot{X} + K(K_1, K_2) Y(X) Y'(X) = 0 \tag{6}$$

with the definition domain written as

$$X \in \{X \mid -L < Y(X) < L\}. \tag{7}$$

## 3. Integrated design method

Based on GSBM, the integrated method of designing nonlinear stiffness systems is proposed, which is an inverse problem in contrast to the constructing process of GSBM discussed above, that is, solving an initial value problem of the first order differential equation about $Y(X)$:

$$\begin{cases} K(K_1, K_2) Y(X) Y'(X) = P(X), \\ Y(0) = \Delta, -L < Y(X) < L, \end{cases} \tag{8}$$

where $P(X)$ is the stiffness force of the system to be designed. For arbitrary stiffness forces, the general solution of Eq. (8) can be given as

$$Y_{1,2}(X) = \pm \sqrt{\Delta^2 + \frac{2}{K(K_1, K_2)} \int_0^X P(X) dX}, X \in \left\{ X \mid -\Delta^2 \leq \frac{2}{K(K_1, K_2)} \int_0^X P(X) dX < L^2 - \Delta^2 \right\}. \tag{9}$$

## 4. Number of mechanical configurations

Geometrically, the mechanical configuration of GSBM is determined by the boundary geometry (roller trajectory) which basically can be divided into two symmetric branches $Y_{1,2}$, as can be seen from Eq. (9). In the case when $\Delta \neq 0$, according to the sign of $K$, $Y_{1,2}$ can be divided into four different forms:

$$\begin{cases} Y_{11} = \sqrt{\Delta^2 + \frac{2}{K(K_1,K_2)}\int_0^X P(X)dX}, X \in \left\{X \mid -\Delta^2 \leq \frac{2}{K(K_1,K_2)}\int_0^X P(X)dX < L^2 - \Delta^2, K(K_1,K_2) > 0\right\} \\ Y_{12} = \sqrt{\Delta^2 + \frac{2}{K(K_1,K_2)}\int_0^X P(X)dX}, X \in \left\{X \mid -\Delta^2 \leq \frac{2}{K(K_1,K_2)}\int_0^X P(X)dX < L^2 - \Delta^2, K(K_1,K_2) < 0\right\} \\ Y_{21} = -\sqrt{\Delta^2 + \frac{2}{K(K_1,K_2)}\int_0^X P(X)dX}, X \in \left\{X \mid -\Delta^2 \leq \frac{2}{K(K_1,K_2)}\int_0^X P(X)dX < L^2 - \Delta^2, K(K_1,K_2) > 0\right\} \\ Y_{22} = -\sqrt{\Delta^2 + \frac{2}{K(K_1,K_2)}\int_0^X P(X)dX}, X \in \left\{X \mid -\Delta^2 \leq \frac{2}{K(K_1,K_2)}\int_0^X P(X)dX < L^2 - \Delta^2, K(K_1,K_2) < 0\right\} \end{cases}$$

(10)

In the case when $\Delta = 0$, according to the sign of $K$, $Y_{1,2}$ can be divided into another four different forms:

$$\begin{cases} Y_{13} = \sqrt{\frac{2}{K(K_1,K_2)}\int_0^X P(X)dX}, X \in \left\{X \mid 0 \leq \frac{2}{K(K_1,K_2)}\int_0^X P(X)dX < L^2, K(K_1,K_2) > 0\right\} \\ Y_{14} = \sqrt{\frac{2}{K(K_1,K_2)}\int_0^X P(X)dX}, X \in \left\{X \mid 0 \leq \frac{2}{K(K_1,K_2)}\int_0^X P(X)dX < L^2, K(K_1,K_2) < 0\right\} \\ Y_{23} = -\sqrt{\frac{2}{K(K_1,K_2)}\int_0^X P(X)dX}, X \in \left\{X \mid 0 \leq \frac{2}{K(K_1,K_2)}\int_0^X P(X)dX < L^2, K(K_1,K_2) > 0\right\} \\ Y_{24} = -\sqrt{\frac{2}{K(K_1,K_2)}\int_0^X P(X)dX}, X \in \left\{X \mid 0 \leq \frac{2}{K(K_1,K_2)}\int_0^X P(X)dX < L^2, K(K_1,K_2) < 0\right\} \end{cases}$$

(11)

For solutions given by Eq. (10), it can be easily proved that they always exist for arbitrary $P(X)$ by intermediate value theorem of continuous function. For solutions given by Eq. (11), under the condition that $\int_0^X P(X)dX < 0$ ($\int_0^X P(X)dX > 0$), $Y_{13}$ and $Y_{23}$ ($Y_{14}$ and $Y_{24}$) do not exist, and under other conditions, $Y_{13}$, $Y_{14}$, $Y_{23}$ and $Y_{24}$ exist simultaneously. Therefore, there are six or eight forms of solutions to Eq. (8), which correspond to six or eight forms of mechanical configurations of GSBM. It is worth noting that GSBM is the most general model of the combination of a linear spring with a nonlinear boundary due to the generality of its boundary (track) and spring (GLSM) forms, such that the mechanical configurations obtained by method Eq. (8) are all of the cases that can be found for the nonlinear stiffness system designed based on the basic combination of a linear spring with a nonlinear boundary.

## 5. An example of softening Duffing system

The general expression of Duffing system with softening nonlinearity can be written as

$$M\ddot{X}+K_3 X^3 =0 (K_3 <0). \tag{12}$$

Combining the stiffness force of system (12) with Eq. (9) yields

$$Y_{1,2}(X) = \pm\sqrt{\Delta^2 + \frac{K_3}{2K(K_1,K_2)}X^4}, X \in \left\{X \mid -\Delta^2 \leq \frac{K_3}{2K(K_1,K_2)}X^4 < L^2 - \Delta^2\right\}, \tag{13}$$

which can be divided into eight valid solutions:

$$\begin{cases} Y_{11} = \sqrt{\Delta^2 + \frac{K_3}{2K(K_1,K_2)}X^4}, X \in \left\{X \mid |X| < \left(-\frac{2K(K_1,K_2)\Delta^2}{K_3}\right)^{\frac{1}{4}}, K(K_1,K_2) > 0\right\} \\ Y_{12} = \sqrt{\Delta^2 + \frac{K_3}{2K(K_1,K_2)}X^4}, X \in \left\{X \mid |X| < \left(\frac{2K(K_1,K_2)(L^2-\Delta^2)}{K_3}\right)^{\frac{1}{4}}, K(K_1,K_2) < 0\right\} \\ Y_{13} = 0, X \in \{X \mid X = 0\} \\ Y_{14} = \sqrt{\frac{K_3}{2K(K_1,K_2)}X^2}, X \in \left\{X \mid |X| < \left(\frac{2K(K_1,K_2)L^2}{K_3}\right)^{\frac{1}{4}}, K(K_1,K_2) < 0\right\} \\ Y_{21} = -\sqrt{\Delta^2 + \frac{K_3}{2K(K_1,K_2)}X^4}, X \in \left\{X \mid |X| < \left(-\frac{2K(K_1,K_2)\Delta^2}{K_3}\right)^{\frac{1}{4}}, K(K_1,K_2) > 0\right\} \\ Y_{22}(X) = -\sqrt{\Delta^2 + \frac{K_3}{2K(K_1,K_2)}X^4}, X \in \left\{X \mid |X| < \left(\frac{2K(K_1,K_2)(L^2-\Delta^2)}{K_3}\right)^{\frac{1}{4}}, K(K_1,K_2) < 0\right\} \\ Y_{23} = 0, X \in \{X \mid X = 0\} \\ Y_{24}(X) = -\sqrt{\frac{K_3}{2K(K_1,K_2)}X^2}, X \in \left\{X \mid |X| < \left(\frac{2K(K_1,K_2)L^2}{K_3}\right)^{\frac{1}{4}}, K(K_1,K_2) < 0\right\} \end{cases}. \tag{14}$$

It can be seen from Eq. (14) that $Y_{13}$ and $Y_{23}$ ($Y_{13}=Y_{23}$) correspond to a single point (0, 0). Under this condition, the corresponding stiffness force is a very special case of the softening stiffness force $P(X) = K_3 X^3 (K_3 < 0)$, i.e., $P = 0, X = 0$. Strictly, the point is the most basic geometric configuration. Without losing generality, $Y_{13}$ and $Y_{23}$ are also taken as a kind of designed boundary form of system (12).

It is assumed that the nonlinear coefficient $K_3$ is set as -5000 Nm$^{-3}$, and the design stiffness $K_d$ and the pre-deformation $\Delta$ are selected as 100 Nm$^{-1}$/-100 Nm$^{-1}$ and ±10 mm/0 mm, respectively. According to the roller trajectory functions given by Eq. (14), the corresponding mechanical configurations of GSBM are plotted in Fig. 7. It can be seen that models (d) and (g) are two special cases of GSBM, i.e. GSBM with $K>0$, $\Delta<0$ and GSBM with $K>0$, $\Delta=0$, which correspond to the mechanical configurations of SSCSM. Benefiting from the general stiffness of GLSM ranging from positive to negative, model (e) ($K<0$, $\Delta<0$) and model (f) ($K<0$, $\Delta=0$) can also be obtained. In addition, the track can be not only compressed but also pulled by GLSM, such that models (a-c) which have symmetric roller trajectories with those of models (d-f) respectively can be obtained as well. For these seven mechanical configurations are all originated from GSBM, they can be seen as homologous septuplets in the physical world reflected with the mathematical model of the designed softening Duffing system.

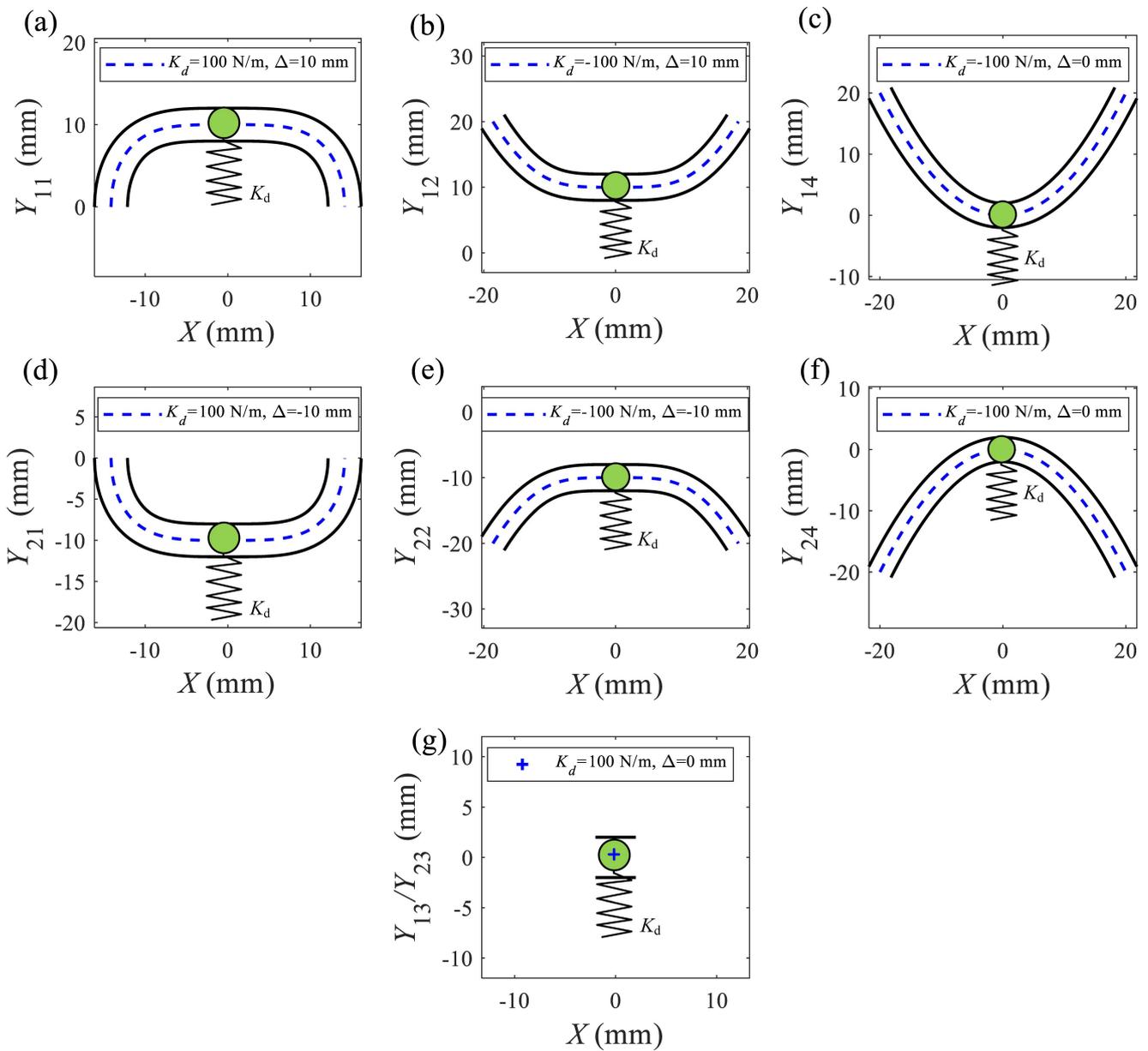

Fig. 7. Seven different mechanical configurations of the designed softening Duffing system.

## 6. Conclusions

In this study, a general spring-boundary model (GSBM) was proposed by introducing a nonlinear track to GLSM. Based on this model, the integrated method to design roller trajectories to customize arbitrary nonlinear stiffness systems was developed. The integrated method considers pre-tensioned, pre-compressed and original length conditions of GLSM with either positive or negative stiffness to design nonlinear stiffness systems. It is proved that the mechanical configurations of nonlinear stiffness systems designed by the combination of a linear spring with a nonlinear boundary are not limited to one, but six or eight forms: for systems with nonnegative or nonpositive potential energy, there are six independent mechanical configurations, and for other systems, there are eight independent mechanical configurations.

## Declaration of Competing Interest

The authors declare that they have no known competing financial interests or personal relationships that could have appeared to influence the work reported in this paper.

## CRediT authorship contribution statement

**Qiangqiang Li**: Methodology, Writing-original draft; **Qingjie Cao**: Writing-review & editing, Project administration, Resources

## Acknowledgment

The authors would like to acknowledge the financial support by the National Natural Science Foundation of China Granted No. 11732006. The authors also acknowledge professor Michel Jean for selflessly providing the English translation of his doctoral dissertation written in 1956.

## References

[1] D. Zou, G. Liu, Z. Rao, et al., A device capable of customizing nonlinear forces for vibration energy harvesting, vibration isolation, and nonlinear energy sink, Mechanical Systems and Signal Processing, 147 (2021) 107101.
[2] D. Ludovico, P. Guardiani, F. Lasagni, et al., Design of Non-Circular Pulleys for Torque Generation: A Convex Optimisation Approach, Ieee Robotics and Automation Letters, 6 (2021) 958-965.
[3] N. Schmit, M. Okada, Design and Realization of a Non-Circular Cable Spool to Synthesize a Nonlinear Rotational Spring, Advanced Robotics, 26 (2012) 235-252.
[4] G. Zhu, Q. Cao, Z. Wang, et al., Road to entire insulation for resonances from a forced mechanical system, Scientific Reports, 12 (2022) 21167.
[5] D. Zou, G. Liu, Z. Rao, et al., Design of a broadband piezoelectric energy harvester with piecewise nonlinearity, Smart Materials and Structures, 30 (2021) 085040.
[6] S. Zuo, D. Wang, Y. Zhang, et al., Design and testing of a parabolic cam-roller quasi-zero-stiffness vibration isolator, International Journal of Mechanical Sciences, 220 (2022) 107146.
[7] M. Jean, Sur les solutions périodique des équations différentielles de la mécaniques, Aix-Marseille Université, 1965.
[8] Q. Cao, Y. Xiong, M. Wiercigroch, A novel model of dipteran flight mechanism, International Journal of dynamics and control, 1 (2013) 1-11.